\documentclass[12pt,reqno]{amsart}
\usepackage{graphicx}%
\usepackage{multirow}%
\usepackage{amsmath,amssymb,amsfonts}%
\usepackage{amsthm}%
\usepackage{mathrsfs}%
\usepackage[title]{appendix}%
\usepackage{xcolor}%
\usepackage{textcomp}%
\usepackage{manyfoot}%
\usepackage{booktabs}%
\usepackage{algorithm}%
\usepackage{algorithmicx}%
\usepackage{algpseudocode}%
\usepackage{listings}%
\usepackage{tikz}
\usepackage{enumitem}
\usepackage{url}

\newcounter{obs}

\newtheorem{theorem}{Theorem}
\newtheorem{corollary}[theorem]{Corollary}
\newtheorem{lemma}[theorem]{Lemma}
\newtheorem{proposition}[theorem]{Proposition}
\theoremstyle{definition}
\newtheorem{definition}[theorem]{Definition}
\theoremstyle{remark}
\newtheorem{remark}[theorem]{Remark}

\numberwithin{equation}{section}
\numberwithin{theorem}{section}

\DeclareMathOperator{\Lip}{Lip}

\begin{document}

\title[Extending graph metrics to hyperspaces of vertex sets]
{Extending graph metrics to hyperspaces of vertex sets: optimization and reinforcement learning applications}

\author[N. Jonard-P\'erez]{Natalia Jonard-P\'erez}
\address{Departamento de Matemáticas, Facultad de Ciencias, Universidad Nacional Autónoma de M\'exico, 04510 Ciudad de M\'exico, M\'exico.}
\email{\textcolor[rgb]{0.00,0.00,0.84}{nat@ciencias.unam.mx}}

\author[E. A.~S\'{a}nchez P\'{e}rez]{Enrique\,A.~S\'{a}nchez P\'{e}rez}
\address{Instituto Universitario de Matem\'atica Pura y Aplicada,
	Universitat Polit\`ecnica de Val\`encia, Camino de Vera s/n, 46022
	Valencia, Spain.}
\email{\textcolor[rgb]{0.00,0.00,0.84}{easancpe@mat.upv.es}}

\author[S. Sanjuan]{Sergi Sanjuan}
\address{Instituto Universitario de Matem\'atica Pura y Aplicada,
	Universitat Polit\`ecnica de Val\`encia, Camino de Vera s/n, 46022
	Valencia, Spain.}
\email{\textcolor[rgb]{0.00,0.00,0.84}{ssansil@upvnet.upv.es}}



\thanks{This research was funded by the Agencia Estatal de Investigaci\'on, grant number PID2022-138342NB-I00.
The research of the first author was funded by PAPIIT grant IN104525 (UNAM, M\'exico). The research of other authors was funded by the European Union's Horizon Europe research and innovation program under the Grant Agreement No. 101059609 (Re-Livestock).}




\begin{abstract}
An undirected connected graph can be interpreted as a metric space by considering the path distance between its vertices. Extending this perspective, we explore how walks in the graph, viewed as sequences of subsets of vertices $V,$ fit into a more general framework where the path distance emerges as a special case of broader metrics that define hypertopologies on the space of subsets of $V.$ In this paper, we propose a general framework to address this issue, motivated by the extension of recent reinforcement learning algorithms on graphs through Lipschitz function extensions. Our approach demonstrates how models for dynamical systems and related machine learning techniques can be adapted to (quasi-pseudo-metric) hyperspaces, offering a novel perspective for the analysis and modeling of dynamical systems. Finally, we illustrate the practical relevance of our theoretical contributions through a concrete graph optimization example and an analytical example involving trajectories of vertex subsets.
 \end{abstract}

\maketitle


\section{Introduction}\label{introduction}

In recent years, the increasing interest in graph theory, driven by its numerous applications across various fields, has led to significant research in graph analysis, often focusing on their treatment as metric spaces. This includes, for example, the general framework for analyzing dynamical systems in trees (connected acyclic undirected graphs) recently introduced by using a nonsymmetric distance function (the bifurcation quasi-pseudo-metric) to model the irreversibility of time-dependent processes by exploiting the asymmetry of these functions \cite{falsan}. 

In this paper, we show that the bifurcation quasi-pseudo-metric is a particular case of a more general and consistent mathematical construction. This approach is based on the idea that any topology on the set of walks in a graph can be viewed as the restriction of a hypertopology on the space of subsets of vertices of the graph, or on the Cartesian product of such sets. Topological models for dynamical systems and machine learning tools can thus be applied within the general framework of quasi-pseudo-metric hypertopologies and associated quasi-pseudo-metrics.

In many learning tasks, however, the basic objects are not single vertices but walks or, more generally, trajectories, that is, ordered sequences of vertex subsets. Viewing trajectories as elements of a hyperspace---and equipping that space with a suitable metric obtained from either the Hausdorff or measure metric---enables the extension of classical Lipschitz-extension techniques to the sequential setting.

Thus, the primary goal of this paper is to extend graph-based metric and learning methods to hyperspaces of vertex subsets and their associated trajectory spaces. In our investigation, we developed a purely set-theoretic approach that enables us to understand distances between subsets of vertices of a graph through an abstract construction based on hypertopologies imposed on sets, defined by classical (quasi-pseudo) metrics on subsets of a metric space (such as the excess function). The study began with the two quasi-metrics on trees introduced in \cite{falsan}, as the specific structure of trees imposes strong limitations on generalizing these metrics to broader classes of graphs. Additionally, \cite{cal} introduced a new metric defined as a modified path distance, resulting in a non-standard metric where points far apart in the path-metric could be close in the new one. However, its computation avoids cycles and requires directed graphs, which is quite restrictive. These limitations further motivated the abstract approach presented here.

\subsection{Graph-Based Models for Dynamical Systems and Hyperspaces}

Machine learning models on graphs often treat the graph as a metric space, with vertices---and sometimes the union of vertices and edges---constituting its elements. The typical metric in a connected graph is defined by the minimum number of steps required to walk from one vertex to another. However, to develop more sophisticated structures and tools, a more abstract framework is necessary. On the other hand, graph theory offers a robust theoretical framework for analyzing discrete dynamical systems and evolution processes, where the primary topological structure is defined by metrics \cite{buck,hak}. A substantial amount of mathematical work has focused on finding appropriate distances; general information can be found in \cite{goy}, while \cite{cao2} addresses metric learning, and reinforcement learning is discussed in \cite{cha,driesetal}. For a more specific approach relevant to this paper, see \cite{kyng} on Lipschitz learning.

When considering the set of finite walks in a graph, each walk can be viewed as an ordered subset of vertices. Therefore, metrics, and more generally topologies, can be regarded as restrictions of hypertopologies defined on the space of all subsets of vertices and edges. 
Motivated by recent developments in reinforcement learning on graphs \cite{falsan}, we adopt this perspective in the present work. While this work considers strategies represented by sequences of states or actions on trees, equipped with bifurcation quasi-metrics or weighted distances, our work extends this approach to trajectories of vertex subsets, using Hausdorff- and measure-based distances and adapting the Lipschitz-extension learning scheme to the resulting trajectory--action spaces.

Thus, in this paper, we mainly analyze general quasi-metric hyperspaces of subsets, and subsequently focus on graph-based metric spaces from both mathematical and machine learning perspectives \cite{entr76}. Specifically, we consider two families of quasi-pseudo-metric hyperspaces---Hausdorff-based and measure-based (Definition \ref{hausdorff_metric} and Proposition \ref{p:measure metric})---and show that they encompass the bifurcation metric on trees (Theorem \ref{t:tree_consistency}). Additionally, we introduce trajectory spaces (Section \ref{trajecmetrics}) and equip them with suitable metrics $\eta_H$ and $\eta_\mu$ (Definitions \ref{hausdorff trajectory} and \ref{measure trajectory}), demonstrating their suitability for sequential decision processes. Moreover, we design a reinforcement learning scheme built upon these metric hyperspaces, relying on Lipschitz extensions of the reward function (Section \ref{sectionRL}). In particular, Section~\ref{sectionRL} identifies sufficient conditions ensuring Lipschitz continuity of the reward on the trajectory--action space and describes how the corresponding Lipschitz extensions are computed in practice.

To illustrate the applicability and utility of our theoretical framework, we present a concrete example of graph optimization (Section \ref{example}). In this example, we explicitly formulate a trajectory-based optimization task on a structured graph as a reinforcement learning problem defined over the introduced quasi-pseudo-metric hyperspaces. 

Moreover, our broader goal is to provide formal applied tools, recognizing that many concepts in this area of research originated from applications in other sciences, such as chemistry \cite{bu14,cheb11,kle93}, as well as formal disciplines like network analysis \cite{bar83,brandes,god11}, routing problems \cite{bell}, and artificial intelligence \cite{gao}.

Finally, let us remark that a key consideration in modeling dynamic problems involving walks in graphs---or any event representable as an (ordered) subset or function within a graph---is their evolving nature. The Markov Decision Processes (MDP) framework \cite{puter} is standard for addressing such issues, and the probabilistic concepts derived from MDPs have been widely applied to graph-based models (see, for example, \cite{puter}). Although this work implicitly utilizes probabilistic ideas from the MDP framework, the methods and arguments presented here are fundamentally different.

\subsection{Theoretical background} 
In general, the elements involved in the evolution processes on graphs are not the graph components themselves -vertices, edges,...- but ordered subsets of these components, which can be associated to actions, trajectories, walks, strategies on the graph, and other characteristic tools that articulate the models. This idea is the cornerstone for extending the framework from graph dynamics to metric hyperspaces of subsets. Similar to systems modeled on graphs, the key tool for constructing machine learning algorithms in this broader setting is the definition of real-valued functions that evaluate ordered subsets (vertices or edges in the case of graphs). Once some values are known for such a function--that represents an evaluation function, a reward function, a characteristic function of a subset---it can be extended using some classical results on extension of real Lipschitz functions on metric spaces after adapting them to the case of metric graphs. 
The use of Lipschitz functions and extensions in Machine Learning goes back to the origin of the subject. Many classical and recent applications can be found in the scientific literature. For example, in \cite{gott,lux} the reader can find some application of Lipschitz extensions in this setting. Very often, a linear space representation supports the mathematical developments, but Lipschitz continuity on general metric spaces---with no linear structure--- are sometimes also considered.
  
According to the systematic use that we make of quasi-metrics,  we do not assume metric symmetry, that is, the distance from an element $x$ to another element $y$ does not necessarily coincide with the distance from $y$ to $x.$ This is particularly convenient for modeling concrete issues of dynamical systems, for example time irreversibility. Some recent mathematical investigations concerning Lipschitz-type functions in this setting are \cite{mustata,mustata2,romsan}, which we use as main references on the topic, together with the book \cite{cob2}.

Let us begin by recalling some concepts on general metric-type definitions. We will denote by $\mathbb N$ the set of natural numbers (containing $0$) and by $\mathbb R^+$ the set of non-negative real numbers. If $\Omega$ is a set, a \textit{quasi-pseudo-metric} in $\Omega$ is a function $d: \Omega \times \Omega \to \mathbb R^+$ such that for all $a, b, c \in \Omega$ the following two properties hold
\begin{enumerate}[label=\roman*.]
	\item \label{cond:eq}$d(a, a) = 0$,
	\item \label{cond:tri} $d(a, b) \leq d(a, c) + d(c, b)$.
\end{enumerate}
In this case, the pair $(\Omega,d)$ is called a \textit{quasi-pseudo-metric space} \cite{kim}.
If $d(a, b) = +\infty$ is allowed, then $d$ is said to be an \textit{extended quasi-pseudo-metric.} If $d$ satisfies conditions \ref{cond:eq}, \ref{cond:tri} and is symmetric, that is
\begin{enumerate}[label=\roman*., resume]
	\item $d(a, b) = d(b, a)$, for all $a, b \in \Omega$,
\end{enumerate}
then $d$ is called a \textit{pseudo-metric}. And, if $d$ satisfies \ref{cond:eq}, \ref{cond:tri} and
\begin{enumerate}[label=\roman*., resume]
    \item $d(a,b)=d(b,a)=0 \implies a=b$,
\end{enumerate}
it is a \textit{quasi-metric} (in the $T_0$ sense). Finally, if $d$ satisfies all of the above conditions then it is said to be a \textit{metric.} 
Depending on the properties satisfied by the function $d$, the pair  $(\Omega, d)$ is called a \textit{pseudo-metric space} a \textit{quasi-metric space}, or \textit{a metric space}, respectively.

In any quasi-pseudo-metric $(\Omega,d)$, it is always possible to define a topology on $\Omega$ using the local neighborhood systems generated by the balls.  If $\varepsilon >0$, we define the \textit{ball of radius} $\varepsilon$ \textit{and center}  $a \in \Omega$ as
$$
B_\varepsilon(a):= \Big\{ b \in \Omega: d(a,b) < \varepsilon \Big\}.
$$
The function $d^{-1}$ defined on $\Omega \times \Omega$ by $d^{-1}(a, b) := d(b, a)$ for all $a,b \in \Omega$ is also a quasi-(pseudo-)metric, called the \textit{conjugate quasi-(pseudo-)metric} of $d$. The function $d^s$ defined on $\Omega \times \Omega$ by $$d^s(a, b) = \max\{d(a, b), d^{-1}(a, b)\}$$ is a pseudo-metric, and a metric if $d$ is a quasi-metric.

A real function $f$ defined on a metric space $(\Omega,d)$ is said to be \textit{Lipschitz} if there is a constant $K>0$ such that $|f(a)-f(b)| \le K \, d(a,b).$ The infimum of such constants $K$ is called the \textit{Lipschitz seminorm} $\Lip(f)$ of $f.$ If $f$ is only defined on a non-empty subspace $S\subset \Omega,$ several classical results allow one to extend $f$ to $\Omega$  preserving the Lipschitz seminorm. In particular, the McShane-Whitney theorem \cite{mcshane,whitney} states that any of the following formulae 
$$
f^M(a):=\sup_{b\in S}\{f(b)-\Lip(f) \, d(a,b)\} \ \ \text{and} \ \ f^W(a):=\inf_{b\in S}\{f(b)+\Lip(f) \, d(a,b)\}, \quad a \in \Omega,
$$
which are known as the \textit{McShane and Whitney extensions,} provide such an extension preserving the constant $\Lip(f).$

In this paper, we will deal with the notion of (quasi-pseudo)-metric hyperspace. Given a topological space $T$, a hyperspace associated to it is a class of subsets of $T,$ that is often required to satisfy certain topological requirements (for example, to be closed in $T$) \cite{cao,kun,kunro}. The study of hypertopologies (topologies defined on hyperspaces) is a well established research area. Metric hyperspaces (i.e., hyperspaces endowed with a metric structure) define a particular class of such general definition, and have been used in many scientific contexts; in particular, they play a relevant role in computer science and artificial intelligence \cite{fuji,pou}. In Section~\ref{met}, we introduce the necessary definitions related to metrics and other distance-like structures on hyperspaces.

The necessary ideas and results on graph theory can be found in general books on the subject. The basic distance commonly used for undirected and connected graphs is the so-called \textit{path distance} (or shortest-path distance), which is calculated as the minimum number of steps needed to reach a vertex (node) from another vertex (see Section 2.2.2 in \cite{brandes}). This distance is a metric if the graph is undirected but becomes an extended-quasi-metric if the graph is directed (in this case, it is necessary to define $d(a,b)=\infty$ if there is no directed walk from $a$ to $b$).
There are various notions of distances between vertices in a graph, and most of them are modifications of the path distance, e.g. by including arc weights in the expression of the metric (see p.297 in \cite{brandes}). The interested reader can find a brief but fairly complete description of the subject in \cite[\S 15]{deza}.
However, there are also other distances worth considering. 

The idea of endowing a graph with a metric dates back to the early development of graph theory.  Some papers in pure mathematics introduced metric notions in the second half of the last century, including the distance between paths (see, for example, \cite{entr76,hak,klein02}). The motivation for introducing such metrics has often come from disciplines other than mathematics. Sociology and chemistry have been the context in which some of the \cite{har53,bar83} graph distances have appeared.
The use of weighted metrics in these contexts is also a topic of current interest (see \cite{che11,god11}). Different metrics for graphs have also been used in the applied framework, as is the case of the resistance distance \cite{cheb11} also related to specific problems in crystallography and social network analysis \cite{bozzo13,kle93}.

\subsection{Quasi-pseudo-metrics for machine learning on graphs}\label{Quasi-pseudo-metrics for machine learning on graphs}

As we said above, the basic graph structure that was studied in the preliminary work \cite{falsan} and what is the starting point of the results of this paper is the tree: this is a connected undirected graph with no cycles. In this case, any couple of vertices can be connected by a unique simple path, and so it was possible to define the bifurcation quasi-metric.  As we will see, we have to consider also the set of ``actions" on the graph for which this metric is provided; in the referred case, these actions are functions acting on walks starting at the first vertex of the tree, which is endowed with the bifurcation metric. It is defined as follows.

Let $G=(V,E)$ be a tree-like graph. A walk is a finite sequence of vertices $u=(u_1,\ldots,u_m)$, with $m\geq1$, such that $u_i\in V$ for all $i$ and $(u_i,u_{i+1})\in E$ for every $i\in\{1,\ldots,m-1\}$. We denote by $|u|=m$ the number of vertices in the sequence, counted with repetitions.
Given two walks $u=(u_1,\ldots,u_m)$ and $v=(v_1,\ldots,v_r)$, we define their longest common initial segment as
$$
u\wedge v=(u_1,\ldots,u_n),
$$
where $n\leq\min\{m,r\}$ is the largest non-negative integer such that $u_i=v_i$ for every $i\in\{1,\ldots,n\}$. If $n=0$, we set $u\wedge v=\emptyset$ and $|u\wedge v|=0$.

\begin{definition}\label{bifurcation}
The \textit{bifurcation quasi-metric} is defined as
$$
q(u,v)=|v|-|u\wedge v|.
$$
The so-called conjugate quasi-metric is then defined as
$$
q^{-1}(u,v)=|u|-|u\wedge v|.
$$
what allows to define a metric as the maximum of $q$ and $q^{-1}$ as
$$
\rho_0(u,v)=\max\{q(u,v),q^{-1}(u,v)\}=\max\{|u|,|v|\}-|u\wedge v|.
$$
\end{definition}
It is proved in Lemma 2 and Proposition 1 of \cite{falsan} that $q$ is, indeed, a quasi-metric.

It should be noted that this bifurcation quasi-metric is defined on the set of possible walks in the graph $G$, rather than on the set of vertices. Therefore, it constitutes a completely different metric structure from the usual one considered on the set of vertices (namely, the minimum path distance described above).

\section{Metric hyperspaces of a quasi-pseudo-metric space} \label{met}

In this section we present two (essentially known) methods for defining metrics on hyperspaces of quasi-metric spaces. Although both of them start with a (pseudo-)quasi-metric on the original space $\Omega,$ the final construction provides metrics on some concrete hyperspaces of subsets belonging to $\mathcal P(\Omega),$ the class of all subsets of $\Omega.$ Let us start with a metric inspired in set theoretic arguments. The second one will be defined using measure theoretic arguments.

\subsection{The Hausdorff quasi-metric on the hyperspace of closed sets.}\label{S:Hausdorff}

Fix a quasi-pseudo-metric space $(\Omega,d),$ and let us write $\mathcal P_0(\Omega)=\mathcal P(\Omega) \setminus \{\emptyset\}$ for the class of all non-empty subsets of $\Omega.$
Consider a subset $A \in \mathcal P_0(\Omega)$. 
For each $x \in \Omega,$ the ``distance" from $x$ to $A,$ denoted by $d(x, A),$ is given by
$$
d(x, A) = \inf\big\{d(x, a) : \, a \in A \big\}. 
$$
The symmetric notion $d(A,x)$ is defined in the obvious way, $d(A,x)=  \inf \big\{d(a, x) : \, a \in A \big\}.$
Note that we can also define a similar function using the conjugate quasi-pseudo-metric, obtaining
$$
d^{-1}(x,A)=\inf\big\{d^{-1}(x,a) : \, a \in A \big\} = \inf\big\{d(a, x) : \, a \in A \big\} , \quad x \in \Omega,
$$
from which we directly get $d^{-1}(x,A)=d(A,x).$

These definitions open the door to the standard construction of a quasi-pseudo-metric on the hyperspace associated with $(\Omega,d),$ and it also allows to extend the notion of metric space to the more general structure called approach space mentioned in the Introduction.

Now we complete the definition of a quasi-pseudo-metric  on the hyperspace $\mathcal P_0(\Omega)$ by means of the original one $d.$
The (extended) quasi-pseudo-metric $d_r: \mathcal P_0(\Omega) \times \mathcal P_0(\Omega) \to \mathbb R^+ \cup \{+\infty\}$ is given by
\begin{equation}\label{excess}
d_r(A,B)= \sup \big\{ d(a,B): \, a \in A \big\}.
\end{equation}
When $d$ is a metric, the value $d_r(A,B)$ is also called the \textit{excess of $A$ over $B$ with respect to $d$} (see, e.g., \cite{Beer}).
In the same way, we define the second (extended) quasi-pseudo-metric
\begin{equation}\label{excess2}
d_l(A,B)= \sup \big\{ d(A,b): \, b \in B \big\} =  \sup \big\{ d^{-1}(b,A): \, b \in B \big\}.
\end{equation}
 
These functions are neither pseudo-metrics nor quasi-metrics, in general. But under some assumptions, the following functional is a metric. 

\begin{definition}\label{hausdorff_metric}
For $A, B \in \mathcal P_0(\Omega),$ we define the extended quasi-pseudo-metric $d_H(A, B)$ as
\begin{equation}
\begin{aligned}
d_H(A, B) &= \max \big\{ d_r(A,B), d_l(A,B) \big\} \\
&= \max \left\{
\sup_{a \in A} \inf_{b \in B} d(a, b), \,
\sup_{b \in B} \inf_{a \in A} d(a, b)
\right\}.
\end{aligned}
\end{equation} 
In general, $d_H$ is an extended quasi-pseudo-metric, which need not be symmetric, since $d_H(\{a\},\{b\})=d(a,b)$. But under some conditions, $d_H$ gives a metric. If $d$ is a metric and $\mathcal C_0(\Omega)$ is the collection of all non-empty closed bounded subsets of $\Omega$, then $d_H$ is a metric on $\mathcal C_0(\Omega).$ In this case, $d_H$ is called the Hausdorff metric. The collection $\mathcal C_0(\Omega)$ topologized by the Hausdorff metric $d_H$ is often called a hyperspace in the context of general topology.
\end{definition}

In this paper, unless stated otherwise, we will work with metrics $d$ and subclasses $\mathcal B(\Omega) \subset \mathcal P_0(\Omega)$ for which $(\mathcal B(\Omega), d_H)$ is a metric space.

\subsection{A metric hyperspace of classes of measurable sets induced by a finite measure: the measure metric.} \label{subsecfin}

Consider a  finite measure space $(\Omega,\Sigma,\mu),$ where  $\Omega $ is a non-empty set and $\Sigma$ is a $\sigma-$algebra. Consider the quotient space $[\Sigma]_\mu$ given by the equivalence classes of measurable subsets that are equal $\mu-$a.e..

\begin{lemma}\label{quasimetric1}
The following formula
$$
q_\mu(A,B)= \mu( B) - \mu(A \cap B), \quad A,B \in \Sigma,
$$
is a quasi-pseudo-metric on $\Sigma,$ and a quasi-metric when the set that is considered is the space of equivalence classes $[\Sigma]_\mu$.
\end{lemma}
\begin{proof}
If $A,B,C \in \Sigma,$ then
$$
\chi_A - \chi_{A \cap C} = \big( \chi_A - \chi_{A \cap B} \big)  + \big( \chi_B - \chi_{C \cap B} \big) + \chi_{A \cap B} - \chi_{A \cap C}- \chi_B + \chi_{C \cap B} 
$$
$$
\le \big( \chi_A - \chi_{A \cap B} \big)  + \big( \chi_B - \chi_{C \cap B} \big).
$$

Note that $q_\mu(A,B)$ can be written as an integral,
$$
q_\mu(A,B)= \int_\Omega \big( \chi_B - \chi_{A \cap B}\big) \, d \mu.
$$
Therefore, integrating the inequality we get
$$
q_\mu(C,A)= \int_\Omega \big( \chi_A - \chi_{A \cap C} \big)  \, d \mu \le   \int_\Omega  \big( \chi_A - \chi_{A \cap B} \big)  d \mu +  \int_\Omega \big( \chi_B - \chi_{C \cap B} \big) d \mu
$$
$$
= q_\mu(C,B) + q_\mu(B,A).
$$

On the other hand, if $q_\mu(A,B)=0=q_\mu(B,A),$ we get that $A \subseteq B$ and $B \subseteq A$ $\mu-$a.e. and so $[A]=[B]$ (where $[A]$ and $[B]$ denote the equivalence classes of $A$ and $B$, respectively,  in $[\Sigma]_\mu$). Therefore, $q_\mu$ is a quasi-metric on $[\Sigma]_\mu$ with conjugate quasi-metric $q_\mu^{-1}([A],[B])= q_\mu([B],[A])$.
\end{proof}

To simplify notation, we will write $q_\mu(A,B)$ and $d_\mu(A,B)$ in place of $q_\mu([A],[B])$ and $d_\mu([A],[B])$ whenever this causes no ambiguity.
The quasi-metric introduced in the previous proposition serves as the basis for defining what we shall refer to as the \textit{measure metric} $d_{\mu},$ given by the following result.

\begin{proposition}\label{p:measure metric}
The function $d_{\mu}: [\Sigma]_\mu \times [\Sigma]_\mu\to \mathbb R^+$ given by
$$
d_{\mu}([A],[B]) = \mu(A \Delta B)
$$
is a metric, called the measure metric.
\end{proposition}
\begin{proof}
Let $[A],[B] \in [\Sigma]_\mu.$
By Lemma \ref{quasimetric1}, it is obvious that the formula 
$$
q_\mu^s([A],[B]) = \max\{q_\mu(A,B),q_\mu(B,A) \}.
$$
gives a metric, that is equivalent (up to constants $1$ and $2$) to the sum $q_{\mu}(A,B)+ q_{\mu}(B,A).$
So, using the notation of the symmetric difference of subsets and choosing the additive version, we get
\begin{equation}\label{measure_metric}
\begin{aligned}
	d_\mu([A],[B]) &= q_\mu(A,B) + q_\mu(B,A) = \mu(B) -\mu(A \cap B) + \mu(A)- \mu (A \cap B) \\
	&= \mu(A \cup B)- \mu(A \cap B) = \mu(A \cup B \setminus A \cap B) = \mu(A \Delta B).
\end{aligned}
\end{equation}
\end{proof}

One interesting detail in comparison with the metric given by the Hausdorff distance is that in this case the (equivalence class of the) measurable set $\emptyset$ can be considered as a set of the metric space $([\Sigma]_\mu, d_{\mu})$ in a natural way.

\section{Hypertopologies on metric graphs}

In order to justify the coherence of our approach, let us prove now that the general framework involving hyperspaces coincides with the reference case of the bifurcation distance on trees introduced in \cite[\S 3.1]{falsan} when we restrict the model to this case.

Let $G=(V,E)$ be a tree-like graph (a tree with a distinguished root $o$) and let $W_0$ be the set of simple paths starting at $o$.
Recall from Definition \ref{bifurcation} the bifurcation quasi-metric $q(u, v) = |v|-|u \wedge v|$ for two walks $u, v \in W_0,$ and the bifurcation metric
\begin{equation*}
	\rho_0(u,v)= \max\{ q(u,v),q^{-1}(u,v)\}= \max\{|u|, |v|\} - |u \wedge v|,
\end{equation*}
where $u \wedge v$ is the common part of the trajectories $v$ and $u$.

Endow the tree-like graph $G=(V,E)$ with the usual path metric $p$ and consider the family $\mathcal S_0=\{S_u:u\in W_0\}$, where $S_u$ denotes the set of vertices of the path $u$. For sets $A, B \in \mathcal S_0$, consider the excess functions defined in Equations \ref{excess}, \ref{excess2} of Section \ref{S:Hausdorff}
\begin{align*}
d_r(A,B) &= \sup \Big\{ p(x,B): x \in A \Big\} \\
d_l(A,B) &= d_r(B,A)
\end{align*}
and the Hausdorff (pseudo-)metric defined in Definition \ref{hausdorff_metric}
$$
d_H(A,B)= \max \Big\{d_r(A,B), d_l(A,B)\Big\},
$$
both induced by $p.$

\begin{theorem}
\label{t:tree_consistency}
Consider a tree-like graph $G=(V,E)$ rooted at $o$ endowed with the usual path metric $p.$ Consider the set of simple paths $W_0$ starting at $o$ and the class $\mathcal S_0$ of all the vertex sets formed by the paths of $W_0.$ Then,
\begin{itemize}   
    \item[i.] the excess functions $d_l$ and $d_r$ are quasi-metrics on $\mathcal S_0$ and coincide, under the identification of each path with its vertex set, with the bifurcation quasi-metrics $q$ and $q^{-1}$ respectively, and

    \item[ii.] $d_H$ is a metric that coincides with $\rho_0.$
\end{itemize}
\end{theorem}
\begin{proof}
    Since $G$ is a tree and all paths in $W_0$ start at the same root, each path is uniquely determined by its vertex set. We now prove the stated identities.
    
    Fix two walks $u=(u_1,u_2,\dots,u_m)$ and $v=(v_1,v_2,\dots,v_r),$ with $u, v \in W_0$ and $u \wedge v = (u_1,u_2,\dots,u_n).$ Then, the bifurcation quasi-metric is defined as
	$$
	q(u, v)= |v| - |u \wedge v| = r-n.
	$$

    Let $S_{u}, S_{v} \in \mathcal S_0$ be the vertex sets formed by the walks $u$ and $v$, respectively. For this particular graph (tree), the excess function $d_l$ is calculated as
	$$
	d_l(S_{u}, S_{v}) = \sup \Big\{ p(S_{u}, v_i): v_i \in v \Big\}.
	$$
	
    Since $p(S_{u}, v_i) = 0$ if $v_i \in S_{u}$ and $p(S_{u}, v_i) = i-n$ if $v_i \notin S_{u}$ (that is, $i>n$), by the uniqueness of paths in a tree, we have that
	$$
	d_l(S_{u}, S_{v}) = |v| - n = r-n = q(u, v).
	$$

	This means that the conjugate also corresponds, $d_r(S_{u}, S_{v}) = q^{-1}(u, v),$ which (by the construction of $d_H$ and $\rho_0$) directly implies that $d_H(S_{u}, S_{v}) = \rho_0(u, v).$
\end{proof}

\section{Metric spaces of trajectories on hyperspaces}\label{trajecmetrics}

Building on the framework established in the previous section, which focused on hypertopologies on metric spaces (particularly on trees), we now generalize these constructions by introducing and defining trajectory spaces as metric hyperspaces, allowing us to extend the theory to more general settings.

Let $\Omega$ be a set and consider its power set $\mathcal P(\Omega).$ Fix a class of subsets $\mathcal B(\Omega) \subset \mathcal P(\Omega)$ and assume that $\mathcal B(\Omega)$ is closed under unions. We now introduce a central concept for our framework, which will serve as a fundamental building block throughout this article.

\begin{definition}
A \textit{trajectory} is a sequence $\varphi=(\varphi_i)_{i=1}^\gamma$, where $\gamma\in\mathbb N\cup\{\infty\}$ and $\varphi_i\in \mathcal B(\Omega)$ for all $i.$
If $\gamma=n\in\mathbb N,$ we say that  $\varphi$ is a \textit{trajectory of $n$ steps.}
We denote by $\mathcal{T}^*(\Omega)$ the set of all (finite and infinite) trajectories, and by $\mathcal{T}(\Omega) \subset \mathcal{T}^*(\Omega)$ the subset consisting of all finite trajectories.
\end{definition}

We are interested in equipping $\mathcal T(\Omega)$ with a distance $\theta$ to construct our target space. In some applications, the set $\mathcal P(\Omega)$ will be substituted by a quotient hyperspace instead, as in the case with the measure-based class $[\Sigma]_\mu$ considered in Section \ref{subsecfin}.

To begin our construction, we need a metric hyperspace $(\mathcal B(\Omega), \rho)$ of subsets of a fixed set $\Omega$, where $\mathcal B(\Omega) \subset \mathcal P(\Omega)$ and $\rho$ is a bounded metric. We are thinking mainly in the cases explained in the previous section. 
As mentioned earlier, $\mathcal B(\Omega)$ is closed under unions and may also be regarded as a quotient of such a hyperspace in specific cases. 

Consider the space $\mathcal T(\Omega) \subset \mathcal T^*(\Omega)$ of all finite trajectories of $\mathcal B(\Omega)$.
For a given trajectory $\varphi = (\varphi_i)_{i=1}^{n},$ we define the $\alpha-$extension function $e_\alpha: \mathcal T(\Omega) \to \mathcal T^*(\Omega)$, with $\alpha \in \mathcal B(\Omega),$ as
$$
e_\alpha(\varphi) = (\varphi_i)_{i=1}^{\infty}, \text{ where } \varphi_i = \alpha \text{ for } i > n.
$$
This function returns a trajectory ``extended to infinity'' by repeating the element $\alpha$.

Fix a padding rule $\alpha(\varphi)\in \mathcal B(\Omega)$ for each finite trajectory, and write $\widetilde{\varphi}=e_{\alpha(\varphi)}(\varphi)$, with $\widetilde{\varphi}=\varphi$ for infinite trajectories. For positive summable weights $W=(w_i)_{i\geq1}$, define on $\mathcal T^*(\Omega)$
\begin{equation}\label{metric_trajectories}
    \theta_{W,\rho}(\varphi^1,\varphi^2):=
    \sum_{i=1}^\infty w_i\,\rho(
    \widetilde{\varphi}^1_i,\widetilde{\varphi}^2_i).
\end{equation}
The series converges by boundedness of $\rho$ and summability. Since symmetry and the triangle inequality hold coordinatewise, $\theta_{W,\rho}$ is a pseudo-metric. With positive weights, zero distance is equivalent to equality of the padded sequences. Let $\mathcal H(\Omega)\subseteq\mathcal T(\Omega)$ be a family of finite trajectories on which the padding map $\varphi\mapsto\widetilde{\varphi}$ is injective. The restriction $\theta$ of $\theta_{W,\rho}$ to $\mathcal H(\Omega)$ therefore defines a metric space $(\mathcal H(\Omega),\theta)$.
To simplify the notation, when there is no risk of confusion about what the sequence $W$ is, we simply denote $\theta_{\rho}$ instead of $\theta_{W, \rho}$.
For finite trajectories, $\rho$ need not be bounded
since the padded tails are constant. Finitely supported
non-negative weights also give a pseudo-metric, which separates
trajectories whenever their positively weighted padded
coordinates determine them uniquely.

\subsection{The space of trajectories with the Hausdorff related metric.}\label{trajectories_Hausdorff}

Let $(\Omega,d)$ be a metric space, and consider a hyperspace $\mathcal C_0(\Omega)$ in which the Hausdorff pseudo-metric $d_H$ defined in \ref{hausdorff_metric} is a metric, that is, the collection of all non-empty closed bounded subsets of $\Omega$.
For the construction of the metric, we propose to use the extension function that repeats the last element of the given trajectory.
\begin{definition}\label{hausdorff trajectory}
Given two finite trajectories $\varphi^1, \varphi^2 \in \mathcal T(\Omega)$ of $n_1$ and $n_2$ steps respectively, and a summable sequence of positive real numbers $(w_i)_i,$ the \textit{Hausdorff pseudo-metric for the space of (finite) trajectories} is defined as
\begin{equation*}
    \eta_H(\varphi^1,\varphi^2)
    :=\theta_{d_H}(e_{\varphi_{n_1}^1}(\varphi^1),
    e_{\varphi_{n_2}^2}(\varphi^2))
    =\sum_{i=1}^\infty w_i\,d_H(
    \widetilde{\varphi}^1_i,\widetilde{\varphi}^2_i),
\end{equation*}
where $e_{\varphi_{n_j}^j}(\varphi^j) = (\varphi^j_1, \varphi^j_2, \dots, \varphi^j_{n_j}, \varphi^j_{n_j}, \dots),$ for $j=1,2$.
\end{definition}

A suitable choice of $\mathcal H(\Omega)$ consists of all trajectories of length $n=1$ or with $\varphi_{n-1}\ne\varphi_n$. The padding is injective on this family: if trajectories of unequal lengths had identical extensions, the longer one would end with two equal elements. Thus $(\mathcal H(\Omega),\eta_H)$ is a metric space.

\subsection{The space of trajectories with the measure related metric.}\label{trajectories_measure}

Let $(\Omega,\Sigma, \mu)$ be a finite measure space. In this case, we consider the class of (equivalence classes of) sets $[\Sigma]_\mu$ and the measure metric $d_\mu$ defined in Equation \ref{measure_metric}.
For the construction of this metric, we propose to use the extension function that repeats the empty set, where $\emptyset$ denotes its measure equivalence class.

\begin{definition}\label{measure trajectory}
Let $\varphi^1,\varphi^2\in\mathcal{T}(\Omega)$ be two finite trajectories of $n_1$ and $n_2$ steps, respectively, and a summable sequence of positive real numbers $W = (w_i)_i.$ Then, we define the \textit{measure pseudo-metric for the space of trajectories} in terms of Equation \ref{metric_trajectories} as
\begin{equation*}
    \eta_\mu(\varphi^1,\varphi^2)
    :=\theta_{d_\mu}(e_{\emptyset}(\varphi^1),
    e_{\emptyset}(\varphi^2))
    =\sum_{i=1}^\infty w_i\,d_\mu(
    \widetilde{\varphi}^1_i,\widetilde{\varphi}^2_i).
\end{equation*}  
\end{definition}

In this case, a suitable choice of $\mathcal H(\Omega)$ consists of all trajectories of length $n=1$ or with $\varphi_n\ne\emptyset$. The padding is injective on this family: identical extensions of unequal lengths would force the longer trajectory to end in the empty-set class. Thus $(\mathcal H(\Omega),\eta_\mu)$ is a metric space.

\begin{remark}
Since we extend all finite trajectories with the empty set $\emptyset$, the number of nonzero terms in the summation defining $\eta_\mu(\varphi^1, \varphi^2)$ is always finite, regardless of the chosen sequence of weights $W = (w_i)_{i}$. Therefore, the weights need only be positive, not summable, for $\eta_\mu$ to be a metric on $\mathcal H(\Omega)$. For example, one can take $w_i=1$ for all $i$.
\end{remark}

\section{Reinforcement learning on hyperspaces: set dynamics}\label{sectionRL}

In this section we develop the reinforcement learning framework associated with the metric space of trajectories $(\mathcal H(\Omega),\theta)$ introduced in the previous section. We first define the corresponding trajectory-action space and establish sufficient conditions ensuring that the reward function is Lipschitz with respect to the natural product metric. This provides a rigorous justification for the use of McShane--Whitney extensions in the learning procedure. We then describe the reinforcement learning scheme based on these extensions.

\subsection{The trajectory-action space and Lipschitz continuity of the reward}

Suppose that we have a metric space of trajectories $(\mathcal H(\Omega),\theta)$ on a metric hyperspace $(\mathcal B(\Omega),\rho),$ where $\theta$ is the restriction to $\mathcal H(\Omega)$ of one of the trajectory distances introduced in the previous section and $\rho$ one of those defined in Section \ref{met}. In this part we define the necessary mathematical elements to construct a canonical reinforcement learning method for dynamics on sets.

Once we have the space $(\mathcal H(\Omega),\theta)$, which represents the set of states and is used to describe the dynamics, we need to define the set $\mathcal A(\Omega)$ of possible actions that would give the next step in the process. That is, if a trajectory $\varphi$ of size $n$ is being developed, the class of (possible) actions are the sets that could be used to define $\varphi_{n+1}.$ In principle, the next step is just either to add another subset of the class $\mathcal B(\Omega)$ to a given trajectory, or to start a new one. The requirements that are imposed for the possible actions depend on the particular process; we require $\mathcal A(\Omega) \subseteq \mathcal B(\Omega).$

In general, the set of possible actions depends on the current trajectory $\varphi,$ and is therefore denoted as $\mathcal A(\Omega)|_{\varphi}.$ In particular, if it depends only on $\varphi_n,$ earlier states do not affect admissibility.
Initial trajectories belong to $\mathcal H(\Omega)$, and every nonterminal admissible action $a$ must satisfy $(\varphi_1,\ldots,\varphi_n,a)\in\mathcal H(\Omega)$. A terminal action ends the episode, and the resulting trajectory is not used for further action selection. In the measure setting, actions, rewards and transitions are defined on measure equivalence classes. We assume finite admissible action sets in the algorithms.
Therefore, in order to formulate the reinforcement learning scheme, the relevant object is not only the metric space of trajectories (states), but the corresponding trajectory-action space
$$
\mathcal G:=\{(\varphi,a)\in \mathcal H(\Omega)\times \mathcal A(\Omega):a\in \mathcal A(\Omega)|_{\varphi}\},
$$
which is naturally endowed with the product metric $\nu((\varphi,a),(\psi,b)):=\theta(\varphi,\psi)+\rho(a,b).$

Recall that a reinforcement learning method is based on the estimate of a reward function $R:$ $\mathcal G$ $\to \mathbb R$ that acts on trajectory-action pairs and is known only on the subset of pairs that have already been experienced during the learning process. The agent then chooses, at step $n$, the action that maximizes the estimated reward $R$ among the possible actions associated with $\varphi$. As will be explained below, and as a generalization of the procedure explained in \cite{falsan}, the method used in this paper is based on extending the reward from the sampled subset of $\mathcal G$ by means of Lipschitz extension formulae of McShane--Whitney type.

For this reason, the key structural requirement is that the reward function be Lipschitz on the metric space $(\mathcal G,\nu)$. This assumption should not be viewed merely as an abstract hypothesis. Rather, our aim is to identify explicit sufficient conditions under which it holds, and to describe classes of graph tasks for which metric continuity of the reward with respect to the trajectory--action space arises naturally; see, for instance, \cite{Hinderer2005,PirottaRestelliBascetta2015,AsadiMisraLittman2018}.

For each finite trajectory $\varphi\in \mathcal H(\Omega)$, we denote by $\widetilde{\varphi}$ the infinite extension used in the definition of $\theta$. Thus, in the Hausdorff setting $\widetilde{\varphi}=e_{\varphi_n}(\varphi)$, while in the measure setting $\widetilde{\varphi}=e_{\emptyset}(\varphi)$.

Let $F:\mathcal H(\Omega)\to\mathbb{R}$ be a trajectory term, let $g:\mathcal A(\Omega)\to\mathbb{R}$ be an action term, and let $h:\mathcal G\to\mathbb{R}$ be a coupling term accounting for the interaction between the current trajectory and the chosen action. We consider reward functions $R:\mathcal G\to\mathbb{R}$ of the form
$$
R(\varphi,a)=F(\varphi)+g(a)+h(\varphi,a).
$$
If $F$ is $L_T$-Lipschitz on $(\mathcal H(\Omega),\theta)$, $g$ is $L_A$-Lipschitz on $\mathcal A(\Omega)$ with respect to $\rho$, and $h$ is $L_C$-Lipschitz on $(\mathcal G,\nu)$, then $R$ is $(L_T+L_A+L_C)$-Lipschitz on $(\mathcal G,\nu)$ by the triangle inequality. Therefore, it is enough to construct Lipschitz trajectory terms, Lipschitz action terms, and Lipschitz coupling terms.

A broad family of trajectory terms is obtained from local scores on the hyperspace $(\mathcal B(\Omega),\rho)$. Indeed, assume that $u:\mathcal B(\Omega)\to\mathbb R$ is $L_u$-Lipschitz and that $(\alpha_i)_{i\ge 1}$ is a sequence of non-negative numbers such that $\sum_{i=1}^\infty \alpha_i < \infty$ and $\alpha_i \le C w_i$ for all $i\ge 1$.
Then
$$
F(\varphi)=\sum_{i=1}^\infty \alpha_i\,u(\widetilde{\varphi}_i)
$$
is well defined for every $\varphi\in \mathcal H(\Omega)$ and defines a $CL_u$-Lipschitz function on $(\mathcal H(\Omega),\theta)$.

We next collect some natural local scores arising in graph problems.

\begin{proposition}
The following assertions hold.

\begin{enumerate}
\item[(i)] Let $D\in \mathcal B(\Omega)$ be fixed. Then
$$
u_D(A):=-\rho(A,D), \qquad A\in \mathcal B(\Omega),
$$
is $1$-Lipschitz with respect to $\rho$.

\item[(ii)] Assume that $(\Omega,d)$ is a metric space and that $\mathcal B(\Omega)\subset \mathcal C_0(\Omega)$ is endowed with the Hausdorff metric $d_H$. If $q:\Omega\to\mathbb R$ is $L_q$-Lipschitz with respect to $d$, then
$$
u_{\sup}(A):=\sup_{x\in A}q(x), \qquad
u_{\inf}(A):=\inf_{x\in A}q(x)
$$
are $L_q$-Lipschitz on $(\mathcal B(\Omega),d_H)$.

\item[(iii)] Assume that $(\Omega,\Sigma,\mu)$ is a finite measure space and that $\mathcal B(\Omega)=[\Sigma]_{\mu}$ is endowed with the measure metric $d_{\mu}$. If $q\in L^{\infty}(\mu)$, then
$$
u_q([A]):=\int_A q\,d\mu
$$
is $\|q\|_{\infty}$-Lipschitz on $([\Sigma]_{\mu},d_{\mu})$.
\end{enumerate}
\end{proposition}
\begin{proof}
\begin{enumerate}
\item[(i)] This follows directly from the reverse triangle inequality for the metric $\rho$.

\item[(ii)] Let $A,B\in \mathcal B(\Omega)$ and let $\varepsilon>0$. For every $x\in A$, there exists $y_\varepsilon\in B$ such that
$$
d(x,y_\varepsilon)\le d_H(A,B)+\varepsilon.
$$
Hence, by the Lipschitz continuity of $q$,
$$
q(x)\le q(y_\varepsilon)+L_q\,d(x,y_\varepsilon)
\le \sup_{z\in B} q(z)+L_q\,d_H(A,B)+L_q\varepsilon.
$$
Taking supremum over $x\in A$, we get
$$
u_{\sup}(A)\le u_{\sup}(B)+L_q\,d_H(A,B)+L_q\varepsilon.
$$
Since $\varepsilon>0$ is arbitrary, it follows that
$$
u_{\sup}(A)\le u_{\sup}(B)+L_q\,d_H(A,B).
$$
Exchanging the roles of $A$ and $B$ yields
$$
|u_{\sup}(A)-u_{\sup}(B)|\le L_q\,d_H(A,B).
$$
The proof for $u_{\inf}$ is analogous.

\item[(iii)] We have
$$
|u_q([A])-u_q([B])|
=\left|\int_{\Omega} q(\chi_A-\chi_B)\,d\mu\right|
\leq \|q\|_{\infty}\int_{\Omega} |\chi_A-\chi_B|\,d\mu.
$$
Since
$$
\int_{\Omega} |\chi_A-\chi_B|\,d\mu=\mu(A\Delta B)=d_{\mu}([A],[B]),
$$
it follows that
$$
|u_q([A])-u_q([B])|
\leq \|q\|_{\infty}\,d_{\mu}([A],[B]).
$$
\end{enumerate}
\end{proof}

Combining the previous observations, we obtain the following family of Lipschitz rewards on the trajectory--action space.

\begin{corollary}\label{rep_reward}
Assume that $u:\mathcal B(\Omega)\to\mathbb R$ is one of the local scores described above, that $(\alpha_i)_{i\geq 1}$ is a non-negative summable sequence satisfying $\alpha_i\leq Cw_i$ for all $i$, that $g:\mathcal A(\Omega)\to\mathbb R$ is $L_A$-Lipschitz with respect to $\rho$, and that $h:\mathcal G\to\mathbb R$ is $L_C$-Lipschitz with respect to $\nu$. Then the reward function
$$
R(\varphi,a)=\sum_{i=1}^{\infty}\alpha_i\,u(\widetilde{\varphi}_i)+g(a)+h(\varphi,a)
$$
is $(CL_u+L_A+L_C)$-Lipschitz on $(\mathcal G,\nu)$.
\end{corollary}

This corollary covers several classes of graph tasks that are relevant in practice. In the singleton case, it yields vertex-based optimization rewards by taking $u(\{v\})=f(v)$ for a Lipschitz function $f:V\to\mathbb R$ with respect to the path metric. It also includes target-reaching and navigation rewards of the form $u(A)=-\rho(A,D)$, where $D$ is a prescribed target set, as well as coverage or exploration rewards in the measure setting, for instance $u([A])=\mu(A\cap D)=\int_A \chi_D\,d\mu$.

The corollary also allows coupling terms. For example, if $P:\mathcal H(\Omega) \to \mathcal B(\Omega)$ is $L_P$-Lipschitz with respect to $\theta$, then
$
h(\varphi,a):=-\rho(a,P(\varphi))
$
is $(1+L_P)$-Lipschitz on $(\mathcal G,\nu)$. This type of term penalizes actions that are far from a trajectory-dependent reference set.

\begin{remark}
The previous conditions are sufficient, not necessary. They identify a broad family of regular graph tasks for which the Lipschitz assumption is structurally justified. By contrast, discontinuous rewards, such as binary rewards given by exact feasibility, exact hitting conditions, or abrupt combinatorial thresholds, need not be Lipschitz with respect to $\nu$; in such cases, a natural alternative is to replace them by distance-based shaping terms. On the other hand, in finite graph problems with bounded horizon and finitely many possible actions at each step, the space $\mathcal G$ is finite, and therefore every reward function $R:\mathcal G\to\mathbb R$ is automatically Lipschitz.
\end{remark}

\subsection{Reinforcement learning scheme}

Once the trajectory-action metric space $(\mathcal G,\nu)$ has been defined and the structural conditions ensuring the Lipschitz continuity of the reward have been identified, the reinforcement learning procedure can be formulated. The method is based on the progressive estimation of the reward function from the sampled data stored in the replay buffer and on the use of McShane--Whitney extensions to evaluate candidate actions.

The proposed scheme estimates immediate rewards for trajectory--action pairs rather than cumulative returns or action-value functions. Action selection is therefore based on a one-step reward estimate, combined with exploration, and does not guarantee global optimality in general.

We now outline the consecutive steps of our reinforcement learning procedure.

\begin{itemize}
\item[(1)] The process starts with an initial state $\varphi_1 \in \mathcal B(\Omega)$ as the starting point of the trajectory $\varphi.$

\item[(2)] An action $a_1$ is chosen uniformly at random from the admissible actions and then performed on the trajectory $\varphi^1$ of size $1,$ given by $\varphi^1_1=\varphi_1.$ The knowledge repository over the system (\textit{replay buffer}) $\mathcal Z(\Omega) \subset \mathcal H(\Omega)\times \mathcal A(\Omega) \times \mathbb R$ is initialized with the triple $(\varphi^1,a_1,R(\varphi^1,a_1)),$ that is
$$
\mathcal Z(\Omega) = \big\{ \big(\varphi^1,a_1,R(\varphi^1,a_1)\big) \big\}.
$$
Recall that this value of the reward function $R(\varphi^1,a_1)$ is already known once the action $a_1$ is executed on the trajectory $\varphi^1.$

\item[(3)] The next state $\varphi_2 \in \mathcal B(\Omega)$ is attained as a result of the application of the action $a_1$ on the trajectory $\varphi^1.$ This gives the next step of the trajectory $\varphi,$ and allows to define the trajectory $\varphi^2$ of size $2,$ given by $\varphi^2_1=\varphi_1,$ $\varphi^2_2=\varphi_2.$

\item[(4)] The action appropriate to $\varphi^2$ is chosen by comparison with the information accumulated in $\mathcal Z (\Omega).$ The algorithm chooses the action $a_2$ for which the reward function estimate is highest, that is, $$a_2 \in \arg\max_{a' \in \mathcal{A}(\Omega)|_{\varphi^2}} \widehat{R}(\varphi^2, a'),$$ where $\widehat{R}$ is the extended reward function. If several actions maximize the estimated reward, ties are broken uniformly at random. Alternatively, an action can be chosen randomly to avoid overlearning.

\item[(5)] In order to choose $a_2$, the estimate $\widehat R$ is computed from the information currently available in the replay buffer. At this stage, the sampled subspace of $\mathcal G$ is
$$
S_1=\{(\varphi^1,a_1)\},
$$
and the reward is known on this set through the value $R(\varphi^1,a_1)$. Therefore, the quantity used by the algorithm is the Lipschitz constant of the restriction of $R$ to $S_1$, namely
$$
\widehat L_1:=\Lip(R|_{S_1})=0.
$$
The estimate $\widehat R$ is then obtained by applying a Lipschitz extension formula to $R|_{S_1}$, for instance through the McShane and Whitney constructions.

\item[(6)] The process continues by adding all the experience of the agent to the replay buffer: once the reward of the pair $(\varphi^2,a_2)$ is observed, the triple $(\varphi^2,a_2,R(\varphi^2,a_2))$ is included in $\mathcal Z(\Omega)$. Thus,
$$
\mathcal Z(\Omega)
=
\big\{
\big(\varphi^1,a_1,R(\varphi^1,a_1)\big),
\big(\varphi^2,a_2,R(\varphi^2,a_2)\big)
\big\}.
$$
Accordingly, the sampled subspace becomes
$$
S_2=\{(\varphi^1,a_1),(\varphi^2,a_2)\}.
$$
At this point, the value of the Lipschitz constant used by the algorithm is updated from the restriction of $R$ to $S_2$. Note that, strictly speaking, each time a new experience is added, the replay buffer should be indexed by the iteration number, that is, written as $\mathcal Z_n(\Omega)$. However, to simplify notation, we keep writing $\mathcal Z(\Omega)$, implicitly assuming that it is updated at each step.

\item[(7)] The procedure follows in this way until a stop criterion (which can be fixed in advance) is satisfied. Then a new trajectory can be initiated, following the same steps as the first one and adding the corresponding experiences to the replay buffer. The metric subspace on which the reward is known is determined by the first two coordinates of the elements stored in $\mathcal Z(\Omega)$.

\item[(8)] This process is repeated by initializing new trajectories and iterating until a global stopping criterion is met, such as a maximum number of epochs (number of trajectories initiated) or convergence to a predefined error threshold. In general, after $n$ records have been stored, the replay buffer determines the sampled subspace
$$
S_n=\{x_i=(\varphi^i,a_i): i=1,\dots,n\}\subset \mathcal G,
$$
and the quantity used by the algorithm is the Lipschitz constant of the restriction of $R$ to $S_n$, namely
$$
\widehat L_n:=\Lip(R|_{S_n})
=
\max_{\substack{1\le i<j\le n\\ \nu(x_i,x_j)>0}}
\frac{|R(x_i)-R(x_j)|}{\nu(x_i,x_j)}.
$$
Observations at zero distance are required to have equal rewards; otherwise, no Lipschitz extension is possible. Such pairs are excluded from the maximum, and we set $\widehat L_n=0$ when no pair has positive distance.

Since the sampled subspaces are nested, $S_{n-1}\subseteq S_n$, the empirical Lipschitz constant can be updated incrementally. For $n\geq2$, we have
$$
\widehat L_n=\max\left\{\widehat L_{n-1},\max_{\substack{1\leq i<n\\ \nu(x_i,x_n)>0}}\frac{|R(x_i)-R(x_n)|}{\nu(x_i,x_n)}\right\},
$$
where the inner maximum is taken to be zero if there are no pairs at positive distance. In practice, the constant is updated before greedy action selection, considering the records added since the previous update. Any remaining records are processed before returning the final predictor in Algorithm~\ref{algorithmRL_train}.

For any candidate pair $x\in \mathcal G$, the McShane and Whitney extensions associated with the sampled data are
$$
R_n^M(x)=\max_{1\le i\le n}\{R(x_i)-\widehat L_n\,\nu(x,x_i)\},
\qquad
R_n^W(x)=\min_{1\le i\le n}\{R(x_i)+\widehat L_n\,\nu(x,x_i)\}.
$$
In practice, the reward estimate is obtained from the convex combination
$$
\widehat R_n(x)=\lambda R_n^M(x)+(1-\lambda)R_n^W(x),
\qquad \lambda\in[0,1].
$$
The empirical Lipschitz constant $\widehat L_n$ satisfies $\widehat L_n\leq \Lip(R)$ whenever $R$ is globally Lipschitz, but need not coincide with it. Consequently, the resulting extensions interpolate the observed rewards but do not necessarily provide bounds for $R$ at unsampled points. The estimate $\widehat R_n$ is therefore used as a heuristic predictor for action selection.

\end{itemize}

In reinforcement learning, if the model always selects the best possible action, it may prematurely assume that certain actions are suboptimal. As a result, it might never explore those actions, even if they are actually the best choices. To mitigate this issue, RL algorithms typically incorporate a random exploration component. For instance, in $\epsilon$ of the cases, the chosen action is selected randomly from all possible actions instead of the one deemed optimal a priori. This exploration rate $\epsilon$ is often defined as a decreasing function over time, starting with a high value to encourage random exploration and gradually decreasing as the model learns.

The training and execution procedures are summarized in Algorithms~\ref{algorithmRL_train} and~\ref{algorithmRL_execute}, respectively, in Appendix~\ref{app:algorithms}. During training, exploration is combined with exploitation to refine the reward estimate, while execution selects the action with the highest estimated reward at each step.

\section[Examples on graphs]{Examples on graphs}\label{example}

\subsection[Vertex-based optimization]{Vertex-based optimization}\label{example_singletons}

Suppose that we have a (finite) connected undirected graph $G = (V, E)$ endowed with the usual path metric $p$, becoming in this way a metric space $(G,p).$ 
Assume that there is a function $f: V \to \mathbb{R}$ and we want to use the proposed learning scheme to search for a vertex $v_{\max}$ for which $f(v_{\max}) = \max\{ f(v) : v \in V \}$.
It is supposed that $f$ has a regular behavior, which in this context would mean that it is Lipschitz with a small Lipschitz seminorm. The elements of the trajectories in this case are the singletons $\{v\}$ defined by the vertices $v \in V$ of the graph.

\begin{figure}[H]
    \centering
    \begin{minipage}[t]{.49\linewidth}
        \centering
        \vspace{0pt}
        \includegraphics[width=\linewidth]{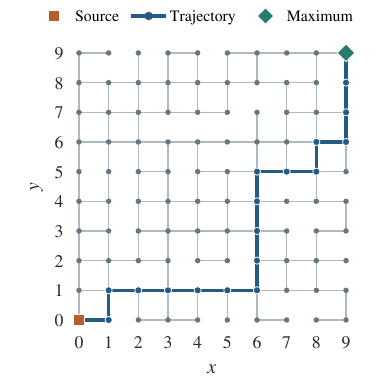}
        \caption{The $10\times10$ grid with $20\%$ of edges removed while preserving connectivity, and the final trajectory.}
        \label{final_trajectory}
    \end{minipage}\hfill
    \begin{minipage}[t]{.49\linewidth}
        \centering
        \vspace{0pt}
        \includegraphics[width=\linewidth]{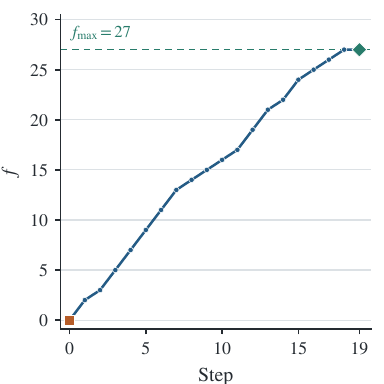}
        \caption{Function values evaluated at each step of the final trajectory.}
        \label{value_f_trajectory}
    \end{minipage}
\end{figure}

As a practical application, we consider the problem of optimizing a real-valued function defined on the vertices of a structured yet partially randomized graph.
The graph $G$ is constructed as a $10 \times 10$ mesh grid in $\mathbb{R}^2,$ with $20\%$ of the edges randomly removed while ensuring connectivity. The graph in Figure~\ref{final_trajectory} can be reproduced using seed $0$ and the graph-generation routine available in the GitHub repository. The coordinates of each node $(x, y)$ are taken as points in $\mathbb{R}^2$.
The function $f$ to be maximized is given by
$$
f(x, y) = 2x + y,
$$
where $(x, y)$ are the coordinates of the corresponding node.
The goal is to find the vertex $v_{\max}$ that maximizes $f$ using the Lipschitz-RL framework. Given the structure of $f$, it follows that the optimal vertex is $v_{\max} = (9,9)$, where the function attains its maximum value of $27$.

The algorithm starts at the source node. The trajectories are constructed by iteratively adding a new vertex to the last one in the trajectory, provided the new vertex is adjacent to the previous one. An action is defined as ``add a new adjacent vertex to the trajectory'', so the set of possible actions, denoted by $\mathcal{A}_{\varphi}$, is determined by the current trajectory $\varphi.$ Specifically, $\mathcal{A}_{\varphi}$ consists of all vertices adjacent to the last vertex in the trajectory that have not been visited, along with the last vertex itself. Selecting the current vertex again is a terminal action. Its reward is evaluated using the current trajectory, and the repeated vertex is not appended. Thus, $\mathcal H(V)$ consists of histories of distinct vertices with length at most $|V|.$ The last vertex of the trajectory generated by Algorithm~\ref{algorithmRL_execute} is taken as the candidate solution.

As a reward function $R$ we can consider a weighted sum of the function $f$ itself over the trajectory. For this example, given a trajectory $\varphi$ of $n$ steps and an action $a\in\mathcal A_\varphi$, we define the reward function by
$$
R(\varphi, a) = \sum_{i=1}^n \frac{2i}{(n+1)(n+2)} f(\varphi_i) + \frac{2}{n+2} f(a).
$$
Here the coefficients depend on $n=|\varphi|$, so the reward is not a literal instance of Corollary \ref{rep_reward}. Nevertheless, $\mathcal G$ is finite and, as shown below, $\nu$ is a metric on it. Hence every reward is Lipschitz on $(\mathcal G,\nu).$

For a fixed trajectory $\varphi$, maximizing the exact reward $R(\varphi,a)$ is equivalent to maximizing $f(a)$. Such a one-step selection may lead to a local maximum. Exploration allows the algorithm to select non-greedy actions and potentially escape local maxima. The proposed procedure is therefore a heuristic without general optimality guarantees.

To extend this function, we use the mean of the McShane and Whitney extensions $$\widehat{R} = \frac{R^M + R^W}{2}.$$

At each training step, the constant $\widehat L$ used in the extension is the empirical Lipschitz constant computed from the current replay buffer with respect to the metric $\nu$ on the sampled pairs $(\varphi^i,a_i)$.

To measure the distance between two trajectories, we use the Hausdorff distance construction explained in Section \ref{trajectories_Hausdorff} with the usual path metric $p.$ For $\varphi^1,\varphi^2\in\mathcal H(V)$, we take $w_i=1/|V|$ for $1\leq i\leq |V|$ and $w_i=0$ thereafter, where $|V|$ is the number of vertices. Thus, the distance between the trajectories is given by
\begin{equation*}
    \eta_H(\varphi^1,\varphi^2) :=\theta_{d_H}(\widetilde{\varphi}^1,\widetilde{\varphi}^2)
    =\sum_{i=1}^{\infty} w_i\,d_H(\widetilde{\varphi}^1_i,\widetilde{\varphi}^2_i)
    = \frac{1}{|V|} \sum_{i=1}^{|V|} d_H(\widetilde{\varphi}^1_i,\widetilde{\varphi}^2_i).
\end{equation*}
Since all histories in $\mathcal H(V)$ have length at most $|V|$, equality of the first $|V|$ padded coordinates implies equality of the entire padded sequences. Terminal padding is injective on these histories by Section~\ref{trajectories_Hausdorff}. Thus the finitely supported weights still give a metric $\eta_H$ on $\mathcal H(V)$ and hence a metric $\nu$ on $\mathcal G$.

For simplicity, we have defined the elements of the trajectory as singletons, meaning that each step $\varphi_i$ consists of a single vertex of the graph. In this formulation, the Hausdorff metric $d_H$ coincides with the usual path metric $p$, as the Hausdorff distance between singletons reduces to the direct distance between the corresponding vertices. In Section~\ref{example_set_valued}, we consider trajectories with non-singleton states.

\begin{figure}[H]
    \centering
    \includegraphics[width=.75\linewidth]{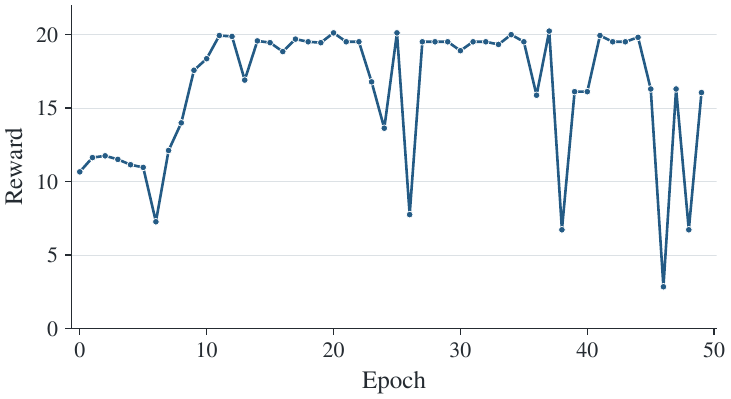}
    \caption{Evolution of the reward value during the training process across $50$ epochs.}
    \label{reward_pregression}
\end{figure}

Using this singleton setup, we implemented the optimization algorithm with the source node $v_0 = (0,0)$ and with an exploration rate of $\epsilon = 0.05$. Figure \ref{reward_pregression} shows the evolution of the reward over $50$ training epochs. After approximately $15$ epochs, the reward frequently reaches values close to $20$, although substantial fluctuations remain, potentially due to exploratory actions under the $\epsilon$-greedy strategy. For the final results, we used $15$ training epochs.

Figure \ref{final_trajectory} depicts the final trajectory generated by the model, showing the path from the source node to the vertex where the function reaches its maximum value. Finally, Figure \ref{value_f_trajectory} presents the function $f(x,y)$ values along the selected trajectory, which reaches the global maximizer $(9,9)$ in this particular run.

\subsection[An example with set-valued states]{An example with set-valued states}\label{example_set_valued}

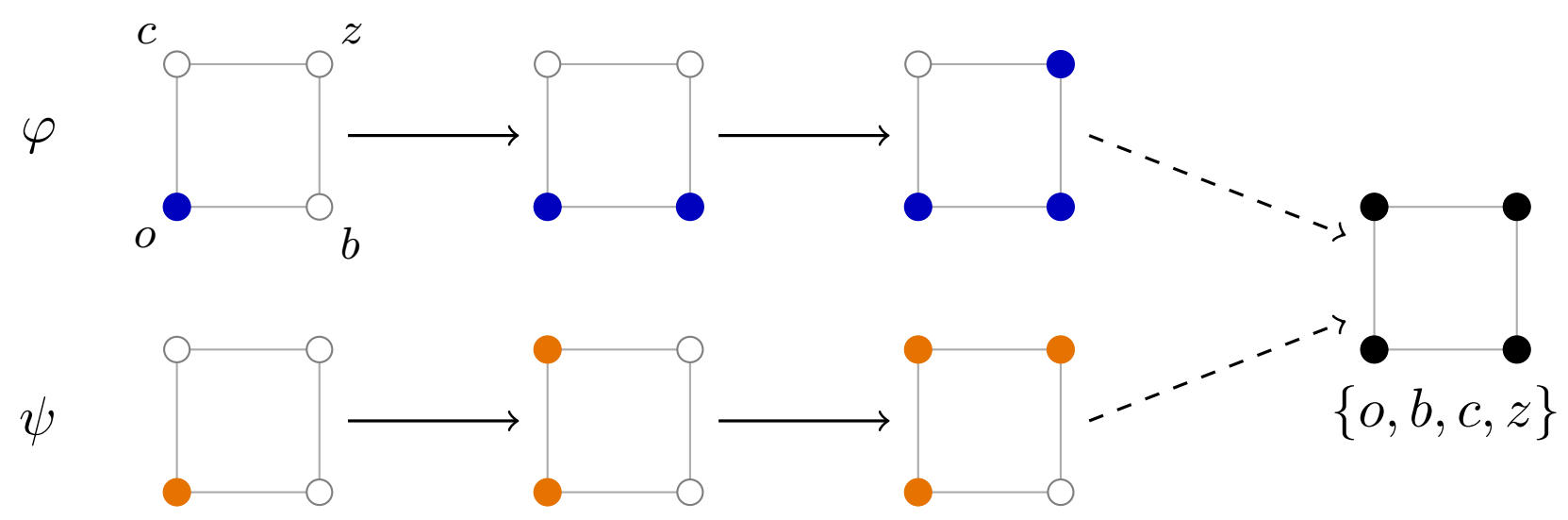
\begin{figure}[htbp]
    \centering
    \begin{tikzpicture}[x=1cm,y=1cm,
        every node/.style={font=\small},
        absent/.style={circle,draw=gray,fill=white,inner sep=1.8pt},
        phiVertex/.style={circle,fill=blue!75!black,inner sep=2pt},
        psiVertex/.style={circle,fill=orange!90!black,inner sep=2pt},
        commonVertex/.style={circle,fill=black,inner sep=2pt}]

        \def\drawset#1#2#3#4{%
            \begin{scope}[shift={(#1,#2)}]
                \draw[gray!65] (0,0) rectangle (1,1);
                \foreach \i/\j in {0/0,1/0,0/1,1/1}
                    \node[absent] at (\i,\j) {};
                \foreach \i/\j in {#4}
                    \node[#3] at (\i,\j) {};
            \end{scope}
        }

        \drawset{0}{2}{phiVertex}{0/0}
        \drawset{2.6}{2}{phiVertex}{0/0,1/0}
        \drawset{5.2}{2}{phiVertex}{0/0,1/0,1/1}

        \drawset{0}{0}{psiVertex}{0/0}
        \drawset{2.6}{0}{psiVertex}{0/0,0/1}
        \drawset{5.2}{0}{psiVertex}{0/0,0/1,1/1}

        \drawset{8.4}{1}{commonVertex}{0/0,1/0,0/1,1/1}

        \node[left] at (-.7,2.5) {$\varphi$};
        \node[left] at (-.7,.5) {$\psi$};

        \node[below left,font=\scriptsize] at (0,2) {$o$};
        \node[below right,font=\scriptsize] at (1,2) {$b$};
        \node[above left,font=\scriptsize] at (0,3) {$c$};
        \node[above right,font=\scriptsize] at (1,3) {$z$};

        \foreach \y in {2.5,.5}{
            \draw[->,semithick] (1.2,\y)--(2.4,\y);
            \draw[->,semithick] (3.8,\y)--(5,\y);
        }

        \draw[->,dashed,semithick] (6.4,2.5)--(8.2,1.8);
        \draw[->,dashed,semithick] (6.4,.5)--(8.2,1.2);

        \node[below=3pt] at (8.9,1) {$\{o,b,c,z\}$};

    \end{tikzpicture}
    \caption{Two trajectories of cumulative vertex sets. Filled vertices belong to the corresponding set. Dashed arrows indicate the transition to a common admissible next state.}
    \label{fig:set_valued_trajectories}
\end{figure}

We now consider trajectories whose elements are cumulative sets of visited vertices rather than singletons. Consider the square graph in Figure~\ref{final_trajectory} with vertices $o=(0,0)$, $b=(1,0)$, $c=(0,1)$ and $z=(1,1)$. The walks $(o,b,z)$ and $(o,c,z)$ define the trajectories
$$
\varphi=(\{o\},\{o,b\},\{o,b,z\}), \qquad \psi=(\{o\},\{o,c\},\{o,c,z\}),
$$
as illustrated in Figure~\ref{fig:set_valued_trajectories}. Although both walks end at the same vertex $z$, their terminal sets are different.

Since $b$ and $c$ are both adjacent to $o$, the vertices belonging to one of the sets but not to the other are at distance one from the other set. Thus, the Hausdorff distances between the corresponding states are
$$
d_H(\varphi_1,\psi_1)=0,\quad
d_H(\varphi_2,\psi_2)=1,\quad
d_H(\varphi_3,\psi_3)=1.
$$
In contrast, the distances between the corresponding vertices of the original walks are $0$, $2$ and $0$, respectively. Applying the terminal extension of Definition~\ref{hausdorff trajectory}, we obtain $\eta_H(\varphi,\psi)=\sum_{i=2}^{\infty}w_i.$ For the corresponding singleton trajectories, the distance would instead be $2w_2$. For instance, taking $w_i=2\cdot 3^{-i}$ gives distances of $1/3$ and $4/9$, respectively. The difference arises because the cumulative states retain information about the previously visited vertices, even when the current vertices coincide.

This distinction can also be incorporated into the reward function. Taking $\varphi_3$ as a fixed target set, the local score $u(A)=-d_H(A,\varphi_3)$ is $1$-Lipschitz by the reverse triangle inequality. Following Corollary~\ref{rep_reward}, we define the reward
$$
R(\chi,a)=\sum_{i=1}^{\infty}w_i\,u(\widetilde{\chi}_i)+u(a).
$$
Since the same weights are used in the reward and in $\eta_H$, this function satisfies
$$
|R(\chi,a)-R(\xi,a')|
\leq \eta_H(\chi,\xi)+d_H(a,a')
=\nu((\chi,a),(\xi,a')).
$$
Both trajectories admit the action $a=\{o,b,c,z\}$, obtained by adding the remaining vertex of the square, as shown in Figure~\ref{fig:set_valued_trajectories}. The initial and second states have scores $-2$ and $-1$, respectively, in both trajectories, while $u(\varphi_3)=0$, $u(\psi_3)=-1$ and $u(a)=-1$. Using the weights defined above, we obtain
$$
R(\varphi,a)=-\frac{23}{9},\qquad
R(\psi,a)=-\frac{8}{3}.
$$
Thus, despite having the same terminal vertex and admitting the same action, the two trajectories receive different rewards according to the sets of vertices visited.

\section{Conclusions}

In this paper, we have proposed a general framework for extending graph-based metric structures to hyperspaces of vertex subsets and to trajectory spaces built from them. This approach provides a unified setting in which several distance constructions on walks and graph dynamics can be interpreted as particular cases of broader quasi-pseudo-metric hyperspace structures. Within this framework, we introduced trajectory spaces endowed with Hausdorff- and measure-related distances and used them to formulate a reinforcement learning procedure based on Lipschitz extensions of the reward function. We also identified sufficient conditions ensuring Lipschitz continuity of the reward on the trajectory--action space and described how the corresponding extensions can be computed in practice.

The graph optimisation example illustrates how the method can guide the search towards an optimal vertex in a structured graph. The numerical experiment should be understood as a proof of concept rather than as a comprehensive empirical validation, and the proposed learning procedure does not provide general guarantees of global optimality. While the numerical experiment is restricted to singleton states, the analytical example illustrates how the Hausdorff metric can be used to compare trajectories of vertex sets and construct Lipschitz reward functions.

\section*{Code availability}

All the algorithms presented in this paper are available in a GitHub repository and can be accessed at the following link: \url{https://github.com/serjj99/Graph_RL.git}.

\begin{appendices}

\section{Pseudocode of the learning procedure}\label{app:algorithms}

\begin{algorithm}
	\caption{Training procedure using Lipschitz reward extension.}
	\label{algorithmRL_train}
	\begin{algorithmic}[1]
		\State \textbf{Set}: Graph $G = (V, E)$
		\State \textbf{Set}: Initialization trajectory $\varphi_1$ 
		\State \textbf{Set}: Number of epochs $N_{\text{epochs}}$
		\State \textbf{Set}: Number of steps per trajectory $N_{\text{steps}}$
		\State \textbf{Set}: Exploration rate $\epsilon$
		\State \textbf{Initialize}: Environment $\varphi \in$ $\mathcal H(V)$
		\State \textbf{Initialize}: Lipschitz extension agent $\widehat{R}$
        \State \textbf{Initialize}: Replay buffer $\mathcal Z(V) \leftarrow \emptyset$, $\widehat L\leftarrow 0$
		\For{each epoch from $1$ to $N_{\text{epochs}}$}
		\State Reset Environment ($\varphi \leftarrow \varphi_1$)
		\For{each step from $1$ to $N_{\text{steps}}$}
		\State Calculate available actions $\mathcal A(V)|_{\varphi}$
		\If{$\mathcal A(V)|_{\varphi} = \emptyset$ (no available actions)}
		\State \textbf{Break step loop}
		\EndIf
        \If{$\mathcal Z(V)=\emptyset$ \textbf{or} Random number $< \epsilon$ (exploration)}
            \State $a \leftarrow \text{Random action of } \mathcal A(V)|_{\varphi}$
        \Else
            \State Update $\widehat L$ incrementally with pending records
            \State Compute $\widehat R$ from the McShane--Whitney formulas using $\widehat L$
            \State $a \leftarrow \arg\max_{a' \in \mathcal{A}(V)|_{\varphi}} \widehat{R}(\varphi, a')$
		\EndIf
        \State Environment: Execute action $a$
        \State $R \leftarrow$ Environment: Get reward $R(\varphi,a)$
        \State $\mathcal Z(V) \leftarrow \mathcal Z(V) \cup \{(\varphi,a,R)\}$
        \If{$a$ is a terminal action}
            \State \textbf{Break step loop}
        \EndIf
        \State $\varphi^+ \leftarrow$ Environment: Get next state
        \State $\varphi \leftarrow \varphi^+$
		\EndFor
		\EndFor
        \State $\widehat L \leftarrow \text{Update incrementally with pending records}$
        \State $\widehat{R} \leftarrow \text{Lipschitz extension of } R \text{ using } \mathcal{Z}(V)$ and $\widehat L$
		\State \textbf{Return} $\widehat{R}$
	\end{algorithmic}
\end{algorithm}

\begin{algorithm}
	\caption{Trajectory execution using the trained model with Lipschitz reward extension.}
	\label{algorithmRL_execute}
	\begin{algorithmic}[1]
		\State \textbf{Set}: Graph $G = (V, E)$
		\State \textbf{Set}: Initialization trajectory $\varphi_1$
		\State \textbf{Set}: Number of steps per trajectory $N_{\text{steps}}$
		\State \textbf{Set}: Lipschitz extension agent $\widehat{R}$ trained by Algorithm \ref{algorithmRL_train}
		\State Reset Environment ($\varphi \leftarrow \varphi_1$)
		\For{each step from $1$ to $N_{\text{steps}}$}
		\State Calculate available actions $\mathcal A(V)|_{\varphi}$
		\If{$\mathcal A(V)|_{\varphi} = \emptyset$ (no available actions)} 
		\State \textbf{Return} $\varphi$
		\EndIf
		\State $a \leftarrow \arg\max_{a' \in \mathcal{A}(V)|_{\varphi}} \widehat{R}(\varphi, a')$
        \State Environment: Execute action $a$
        \If{$a$ is a terminal action}
            \State \textbf{Return} $\varphi$
        \EndIf
        \State $\varphi^+ \leftarrow$ Environment: Get next state
        \State $\varphi \leftarrow \varphi^+$
		\EndFor
		\State \textbf{Return} $\varphi$
	\end{algorithmic}
\end{algorithm}

\end{appendices}

\clearpage

\bibliographystyle{amsplain}
\bibliography{biblio}

\end{document}